\documentclass[12pt]{amsart}
\usepackage{amssymb, latexsym}
\pdfoutput=1
\usepackage{graphicx}
\usepackage{listings}
\begin{document}

\title{Analysis Of A Long Memory Circular Convolution Model}
\author{Robert Kimberk}
\address{Smithsonian Astrophysical Observatory, Cambridge  Massachusetts  }
\email{rkimberk@cfa.harvard.edu}

\begin{abstract}
A  stochastic  model, the product of a circulant matrix and a random normal vector, is shown to produce an evolutive long memory time series with a power law spectral density. The distribution of the time series,  a beta location scale family of distributions, yields a connection to the unit centered spherical distribution and directional statistics. The eigenanalysis of the deterministic circulant matrix is demonstrated to provide estimates of the discrete Fourier spectral trend, the intrinsic dimension, the probability density shape parameter of the resulting time series, the condition number of the matrix, and a  principle component analysis. Examples of the R code, used as the exploratory element of the work, are given as constructive elements of the paper. The R code may be copied,  pasted into a R editor, and explored.
\end{abstract}

\maketitle

\keywords{\bf {Keywords:} long memory time series, circular convolution, circulant matrix, eigenanalysis, intrinsic dimension, direction statistics, beta distribution, unit  centered d dimensional spherical distribution}

\section{Introduction}
Peter M. Robinson \cite{pR03}  stated that it is commonly understood that the discrete  Fourier spectrum of a long memory time series has a pole (is divergent) at frequency zero. Rafal Weron \cite{rW02} equates long memory to $1/f$ noise, meaning the power spectral density has a $|frequency |^{-\beta}$ power law curve, which, for $\beta > 0$,  is also divergent at frequency zero.

This paper adopts the $|frequency |^{-\beta}$ definition of a long memory time series, where $\beta$, the slope of the log log graph of the power spectral density trend, can assume any value between zero and ten. However, as all sets of quantities  described in this paper are bounded in both length and magnitude, the lowest discrete Fourier spectral frequency may approach zero as the length of the model increases,  but never becomes zero. The frequency approach to zero, for time series with $\beta > 0$, produces an increase in the power spectral density,  resulting in  an increase of variance as the length of the time series increases. The long memory time series is non stationary with respect to variance.

This paper's model uses a circular convolution \cite{dP99} to generate the long memory time series. The eigenvectors of the circulant matrix that generates the convolution are the discrete Fourier transform basis vectors and the eigenvalues of the matrix are the real coefficients of the discrete Fourier transform.

 By the end of the paper, it will be shown that both the intrinsic dimension, and the probability density function, of the long memory time series are functions of the number of significant  low frequency spectral components of the Fourier transform of the time series, and the corresponding eigenvalues of the circulant matrix.

\section{ The R Language Model that Produces Long Memory Sample Realizations}

\begin{verbatim}

# long memory time series
# load package named "magic" to produce circulant matrix
beta = 7
n = 500
seed = 5
rn = 0
q = seq(from = -1/2, to = -1/n, by = 1/n) 
r = 1/(2*n)
s = -rev(q)
frequency = c(q,r,s)
rn = length(frequency)
density = (abs (frequency))^(- beta/2)
b = Mod((rn^(-1/2)) * fft(density, inverse =TRUE))
# set.seed( seed)
epsilon = rnorm(rn, mean = 0, sd = 1)
circ = (rn^(-1/2)) * circulant (b)
longmem =  circ  %*%  epsilon
plot(longmem, type="l", lwd =2)
\end{verbatim}

\begin{figure}[h]
\centering \includegraphics [scale = .62]{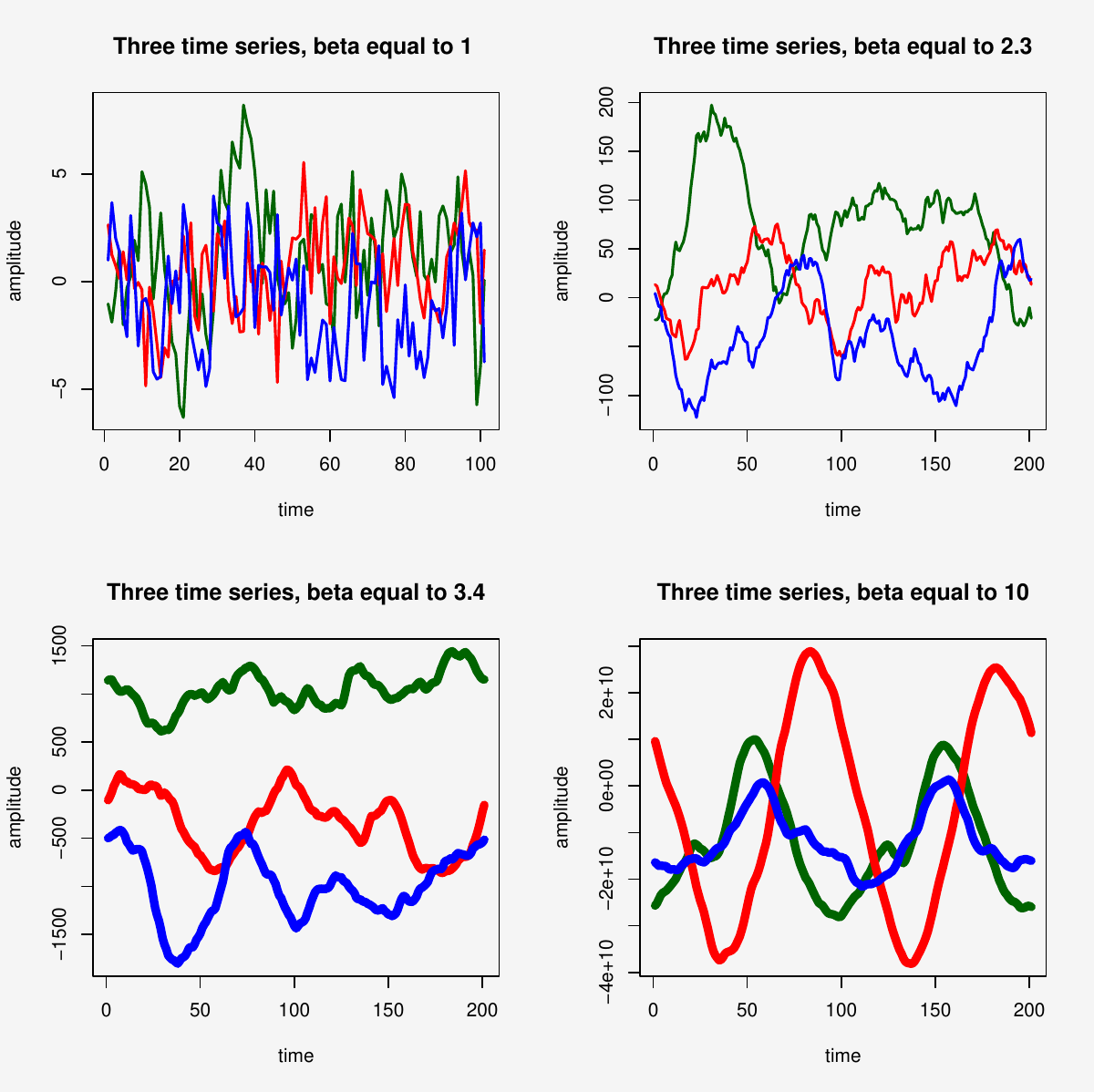}
\caption {Long Memory Time Series From Circulant Matrix Model}
\label {Fi:001}
\end {figure}

\subsection{The Model Components Explained}\     
\\

\noindent Load a required R language package, and variables are set by user:
\begin{itemize}
	\item Load the R language package "magic" to allow production of circulant matrix.
	\item The variable "beta", in the interval zero to ten , the slope of the log log plot of the power spectral density of the time series, is set by user.
	\item The integer variable "n", the reciprocal of the increment of the spectral range of -1/2...-1/n, 1/2n,1/n...1/2 is set by user. This is the Nyquist spectral range of a time series sampled at unit time intervals.
	\item If set.seed( )  is to be used, set seed to an integer value. 
	\item The measured lengths of the spectral range and resulting time series, "rn", is initialized to 0 . If "n" is odd "rn" will be equal to "n". if "n" is even "rn" will be "n"+1. "rn" will be set later  by the program.
\end{itemize}  

\noindent The construction of the variable "frequency" which is the spectral range.
\begin{itemize}
	\item Variable "q" is the vector of negative frequencies.
	\item Variable "r" is the center frequency $1/2n$.
	\item Variable "s" is the vector  of positive frequencies. 
	\item The variable"frequency" is the concatenation of "q","r", and "s",  the Nyquist spectral range.
\end{itemize}

\noindent  Deriving the first row of the circulant matrix:
\begin{itemize}
	\item Variable "rn" is set to the measured length of "frequency".
	\item The variable "density" is set to $|frequency|^{-\beta/2}$. "density" is the sequence of the moduli of  Fourier coefficients of the time series.
	\item The variable  "b" is the moduli of the inverse unitary discrete Fourier transform of "density". 
\end{itemize}

\noindent The normal random vector generated , the circulant matrix produced, and the matrix vector product produces the long memory time series.
\begin{itemize}
	\item Set.seed( ) may be uncommented by removing the hash sign \# to allow multiple extended sequences of "longmem".
	\item The normal random vector "epsilon"$\sim$ N(0,1) is produced.
	\item The circulant matrix "circ" is produced from "b" using the "circulant" command from the package "magic", and normalized by one over the square root of "rn".
	\item The matrix vector product  of "circ" and "epsilon" completes the circular convolution that produces the long memory time series "longmem". A plot of "longmem" is created.
	
\end{itemize}

\section{The Symmetric Circulant Matrix}

The first row of the circulant matrix "circ", produced by the model parameters  "rn" = 7 and "beta" = 3 is: \\

(12.510,    8.253,    6.140,    5.543,    5.543,    6.140,    8.253)
\\

The sequence, after the left most element, is a palindromic sequence of numbers, such as 1,2,3,3,2,1 or 0,0,0,0 or 1,1,1,1. The second row of the circulant matrix  is the first row  circularly shifted to the right, with the last element in the row moved to the first position. The next rows follow the same pattern of a circular shift of the row above. An example of a symmetric circulant matrix of letters is : \\

\begin{equation*}
\left [
\begin{matrix}
		A & B & C & B\\
		B & A & B & C\\
		C & B & A & B\\
		B & C & B & A\\
\end{matrix}
\right ]
\end{equation*}\\

The matrix is symmetrical about its main diagonal.The matrix is  self adjoint,  equal to its transpose. The first row of the matrix is identical to the first column, the second row is identical to the second column, and so on. The square symmetric positive definite  circulant matrix used to model long memory time series in this paper has an odd number of columns and rows. The matrix is multiplied by $rn^{-0.5}$ to produce an unitary transformation.

\subsection{The Top Row of The Circulant Matrix} \
\\
\begin{table}[h]
\begin{center}
\begin{tabular}{ | l |  r | r | }
\hline
"beta"	&  top row "circ"				& 	"density"	\\ \hline
0	&  1,  0,  0, 0, 0				&  	1, 1, 1, 1, 1	\\ \hline
7	&  664, 637, 612, 612, 637		&	11.3,  67.6,  3162,  67.6, 11.3	\\  \hline
\end{tabular} 
\linebreak
\linebreak
\caption{The Response of Fourier Transform Pair "Circ" top row, and "Density" to "Beta"}\label{Ta : first}	
\end{center}
\end{table}

Table 1 is a chart of the top row of the normalized matrix "circ" and the variable "density" in response to two values of "beta". The length of the model is five. The top row of "circ" and the variable"density" are a Fourier transform pair, giving  insight to the response of the model over the range of "beta". 

When "beta" is zero,  the the frequency domain variable  "density" is a vector of numbers raised to the zero power, which is a vector of ones. The inverse Fourier transform of "density" is the time domain  top row of the matrix "circ", a vector of one followed by zeros. The matrix "circ", after the circular permutation of the top row, becomes an identity matrix and the resulting time series after circular convolution,  is the vector "epsilon". 

When "beta" is seven, the frequency domain variable "density" is a vector with large value at nearly zero frequency and smaller values  symmetrically located on either side. The time domain variable, the top row of "circ" is a vector of nearly similar values. The matrix "circ", after the circular permutation of the top row, becomes a nearly singular matrix with large condition number.

\section{Equivalence Class of Location and Scale Transformationed Distributions }\

Michel Loeve \cite{mL55} stated in "Probability Theory" that the type or law of a distribution function is independent of location and scale. Expressed differently, the probability distributions of random variables differing only by location and or scale transformations are an equivalence class. The terms distribution and density in this paper refer to this equivalence class. Some items to remember are:

\noindent
\begin{itemize}
\item The dot product of two vectors A  and B is equal to the product of the  Euclidean $l^2$ norms of A and B and the cosine of the angle between A and B. The dot product of A and B = $ \parallel A \parallel \parallel B \parallel cos \theta$.
\item Each component of vector "longmem" is the dot product of a row of the matrix "circ" and the vector "epsilon". 
\item The norms of each row of "circ" are equal. Each row of "circ" is composed of a permutation of the same numbers.  
\end{itemize}

\begin{equation}
"longmem" = \\ \parallel "circ_{row 1} \parallel \parallel "epsilon" \parallel
\left[
\begin{matrix}
cos \theta_{1}\\
\vdots \\
cos \theta _{n}\\
\end{matrix}
\right]
\end{equation} \

 Dividing the vector "longmem" by the product of the two $l^2$ norms in equation (1) yields an equivalent location scale distribution. 

It is interesting to note that the  probability distribution of the vector of cosines of angles between the unit magnitude rows of "circ" and the unit magnitude vector "epsilon" is a directional statistic. \\

It will be shown in section 6 that the probability distribution of the vector of cosines has a general beta distribution as defined by NIST ( National Institute of Standards and Technology) \cite{NIST}  on their online EDA Handbook ( John Tukey's Exploratory Data Analysis). NIST defines the lower limit of the values of the general beta distribution, "a", as a location parameter, "b" as the upper limit, and the range, "b - a" ,  as the scale parameter. In order to arrive at a standard beta distribution, defined as having support on the interval [0,1], each component of the cosine vector in equation (1)  has location parameter "a" subtracted from it to produce a lower limit of  zero, then the scale parameter "b-a"  divides the result of the subtraction to produce an upper limit of one. The result is shown  in equation (2).

\begin{equation}
\left[
\begin{matrix} 
 (cos\theta_{1} - a)/ (b -a) \\
\vdots \\
(cos \theta_{n} - a)/ (b-a)
\end{matrix}
\right]
\end{equation} \

\section{Histograms of Long Memory Time Series at different values of "beta"}

The vector of equation (2)  is a small sample of an autocorrelated time series, so in the following R program, many standardized vectors of equation (2)  will be concatenated to produce an adequate representation of the probability density function.

\subsection{ The R Language Program to Produce A Histogram}

\begin{verbatim}

#distribution of a long memory time series
# load R package "magic"
beta = 2.3
n = 200
result = 0

q = seq(from = -1/2, to = -1/n, by = 1/n)
r = 1/(2*n)
s = - rev(q)
frequency = c(q,r,s)
rn = length(frequency)
density = (abs(frequency))^(-beta/2)
b = Mod((rn^(-1/2)) * fft(density, inverse = TRUE))
circ = (rn^(-1/2)) * circulant (b)

for(i in 1:200) {
epsilon = rnorm(rn, mean = 0, sd = 1)
norm_epsilon = (sum(epsilon^2))^0.5
norm_row = (sum(circ[1, ]^2))^.5
longmem = circ %*% epsilon
cosvec = longmem/ (norm_epsilon * norm_row)
standard = (cosvec - min(cosvec))/(max(cosvec) - min(cosvec))
result = c(result, standard)  }

result = setdiff(result, c(0,1))
hist(result, freq = FALSE, breaks = 100)

\end{verbatim}

\subsection{The Histogram Program Explained} \

\noindent Much of this code is identical to the code in section 2.1 . The explanation will begin at the for loop staring with "for(i in 1:200)".
\begin{itemize} 
	\item The for loop "for(i in 1:200)" starts a for loop from 1 to 200, indexed by i. For each pass through the loop a new vector "epsilon " is created. The next two line calculate the vector norms of "epsilon"  and the top row of the circulant matrix "circ", these are  called "norm epsilon" and "norm row".
	\item The time series "longmem" is created from the matrix product of "circ" and the vector "epsilon". On the next line of code, the vector of cosines "cosvec"  is created by dividing "longmem" by the product of " norm epsilon" and "norm row". 
	\item  "cosvec" is standardized as in equation 2 to produce the vector "standard".
	\item The 200 different  vectors "standard" produced by the for loop are  concatenated as the  vector "result". The vector "result" then has the excess ones and zeros removed by the setdiff(result, c(1,0) command.
	\item The histogram, an empirical density function, is produced with 50 bins from "result".
 
\end{itemize}

\begin{figure}[h]
\centering \includegraphics [scale = .62]{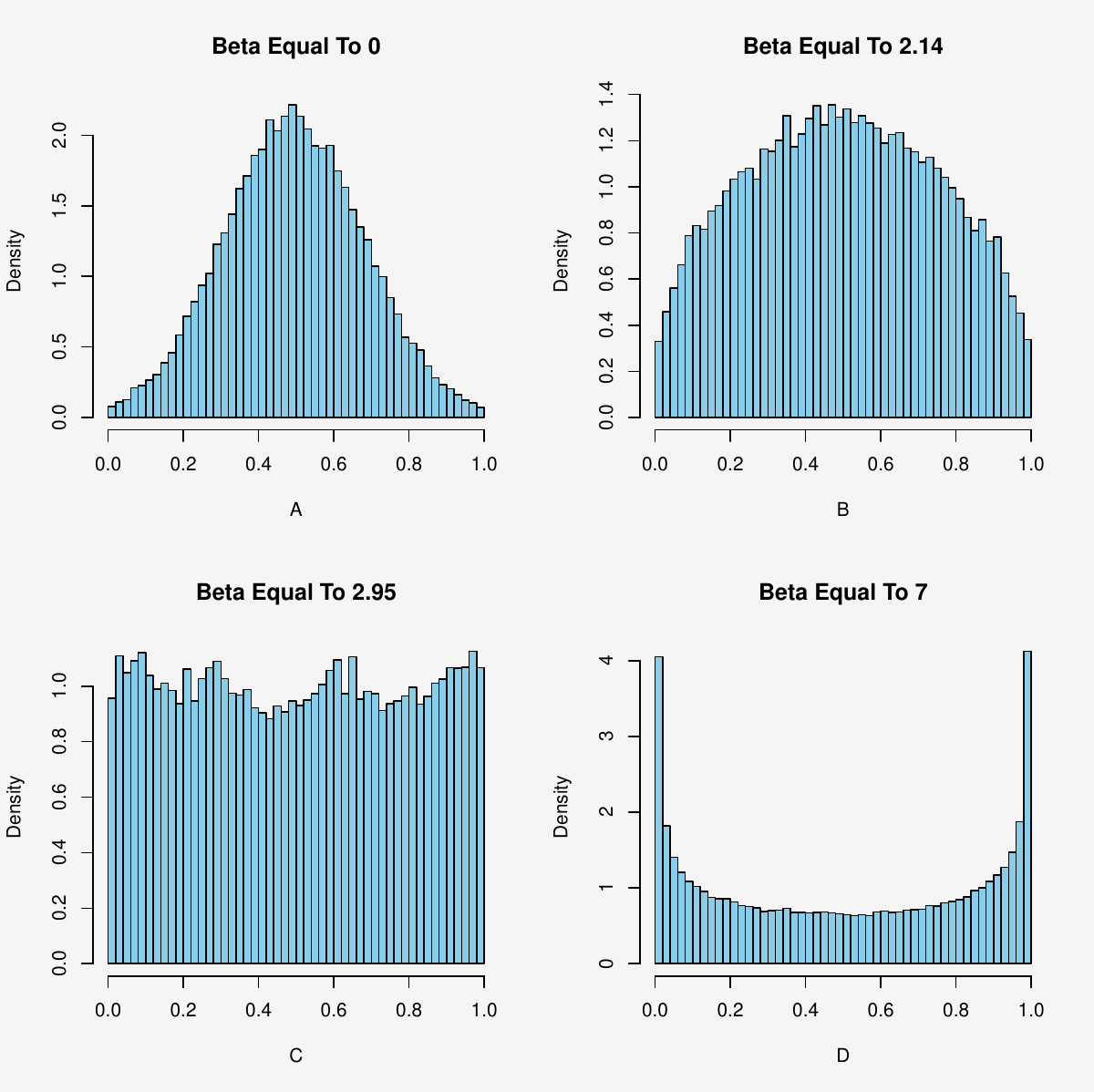}
\caption {Four Histograms Of Time Series}
\label {Fi:002}
\end{figure} 

\section{The General Beta Distribution Of Long Memory Time Series} 

The  histograms of figure 2 are generated from the paper's model of long memory time series with different values of the slope of the power spectral density "beta". They represent  four standardized density functions in the beta family of distributions. The top left is the truncated normal, the top right is the Wigner semicircle, the bottom left is the uniform, and the bottom right is the arcsine distribution. These distributions are not unique to the model used in this paper, the distributions can be generated by the R language ARFIMA model of Granger and Joyeux as well as the R language TK95 model of Timmer and Konig. Details can be found in an earlier paper on arXiv by the author of this paper.

\subsection{The Ratio Of Variance Divided By The Range Squared} \

The variance of the general beta distribution is given by NIST as:

\begin{equation}
\sigma^2  = (b -a )^2 \frac { \alpha \beta}{ ( \alpha +\beta)^2 ( \alpha + \beta +1)}
\end{equation}

Where, a is the lower limit of the range, b is the upper limit of the range and $(a - b)^2$ is the range squared.  The variables, $\alpha$ and $\beta$, are the two shape parameters of the beta distribution. Fortunately for this paper, the beta distributions in this paper are symmetric, and $\alpha = \beta$, so all references to $\beta$ are replaced by $\alpha$. Dividing both sides of equation (3) by the range squared, replacing $\beta$ with $\alpha$ and simplifying yields:

\begin{equation}
\frac {\sigma^2 } {( b - a )^2} = \frac{1}{8 \alpha + 4}
\end{equation}

A property of the left side of equation 4 is scale invariance. Apostal \cite{tA65} states that al  homogeneous function $f$ of degree $p$ is a function of a sequence $x$ and a scalar $t$ such that :

\begin{equation}
f(t(x)) = t^p f(x)
\end{equation}

Variance and the range squared are both homogeneous functions of degree two, so any scale transformation of the the random sequence $x$ cancels in the ratio. Variance and range are also location invariant. The result is that the variance to range squared ratio, easily measured from finite long memory sequences , can be used to determine the beta shape parameter $\alpha$.

\begin{table} 
\begin{center}
\begin{tabular} { |  l  |  r |  r | }
\hline
 Density 		&	Variance / Range Squared				&  Shape Parameter Alpha		\\  \hline
Normal		&	 				 Approaching  0			&	Approaching Infinity 	\\  \hline
Wigner 		&						1/16				&	3/2		\\  \hline
Uniform		&						1/12				&	1		\\  \hline
Arcsine		&						1/8				&	1/2		\\  \hline
Bernoulli		&						1/4				&	0		\\  \hline
\end{tabular}
\linebreak 
\linebreak
\caption{Table of Variance Divided By Range Squared, and Beta Distribution Shape Parameter }\label{ta : first}
\end{center}
\end{table}

\subsection{Popoviciu's Inequality of Variance}

Tiberiu Popoviciu's variance inequality on bounded probability distributions \cite{lG17}, with (b - a) equal to the range of the distribution, states:

\begin{equation}
\sigma^2  \leqq \frac{ (b-a)^2 }{4} 
\end{equation}\

with only the Bernoulli distribution having equality with 1/4. This inequality provides an upper bound of the ratio of $\sigma^2$/$(b-a)^2$ .

A proof of Popoviciu's theorem as applied to symmetric distributions  follows.\

Note that  a symmetric Bernoulli distribution with probability 1/2, has a variance greater than any other bounded distribution with the same range  (b - a). The distance from the extreme ends of the Bernoulli distribution, where the entire mass of the random variable values are found,  to the mean is $x-\mu$, eqaul to plus or minus (b - a)/2 .
The ratio of variance to squared range  of the Bernoulli distribution with probability 1/2 is:

\begin{align}
\sigma^2 &= 1/2\left[ \left( \frac{(b - a)}{2} \right) ^2 + \left (\frac{(a - b)}{2} \right)^2 \right] = \frac{(b - a)^2}{4} \\
\frac{\sigma^2} {range^2} &=  \frac{(b - a)^2}{4} \frac {1}{(b - a)^2} = 1/4  \notag  
\end{align} \\

The  ratio of variance to range squared of  any bounded distribution   is $\leqq 1/4$.\\

\subsection{ The Unit  Centered  Sphere and Cosine Similarity} \

The uniformly distributed unit vector variables "epsilon" in the long memory histogram model  can be represented as a $S^{d-1}$ spherical surface, with d $\in \mathbb{N}$.  The unit row vectors of the circulant matrix are also members of the unit hypersphere.  What determines the particular beta distribution is the intrinsic dimension d of the dot product of the unit  rows of the circulant matrix and the unit normal vectors "epsilon". Four references derive the beta distribution on the unit d dimensional sphere. Djalil Chafai's two papers \cite{dC22} \cite{dC21}  use the Funk Hecke formula to derive the beta density function. William A Huber \cite{wH14} uses a geometric construction on the $S^{d-1}$ sphere to arrive at the same density function. Tony Cai and Tiefeng Jiang \cite{tC11}  derive the density function from the distribution of Pearson correlation coefficients.  The density function given by all four references is  the following expression:

\begin{align}
f_{t} &= \frac{ \Gamma (\frac{d}{2}) }{ \pi^{1/2 } \Gamma(\frac{d-1}{2})}  (1-t^2)^{\frac{d-3}{2}}     \\ \notag
f_{t} &\propto (1-t^2)^{\frac{d-3}{2}}
\end{align}
\
The random variable t in equation (8), created by the many dot product cosines, has a support of [-1,1]. An affine transformation of variable  t to variable u = (t+1)/2 will have the standard beta distribution support of [0,1].  Below is a template of the  linear transformation of the density function of $f_{t}$  to $ f_{u}$. 
\begin{align}
u &= at + b  \\ \notag
t &=  (u - b)/a  \\  \notag
f_{u} &= 1/|a| \ f_{t}\left [ \frac{ u - b}{a} \right ]  \notag
\end{align}

Insert the values of a = b = 1/2 in the second line of equation (9), where  t = 2u - 1  and - t =  1 - 2u. Substituting 1-2u for -t in equation (8) and simplifying the expression  yield a beta distribution with support  on the interval [0,1].

\begin{align}
t &=  (u - b)/a , \ a=b=0.5 \\  \notag
t &= 2u -1 \\ \notag
-t &= 1 -2u \\ \notag
f_{t} &\propto (1-t^2)^{\frac{d-3}{2}} \\ \notag
f_{u}&\propto 2 \left [ 1 - ( 1-2u)^{2} \right ]^{\frac{d-3}{2}} \\ \notag
f_{u} &\propto 2 \left [ 1- (4u^{2}- 4u +1) \right ]^{\frac{d-3}{2}} \\ \notag
f_{u} &\propto  2 \left [ 4 ( u - u^{2}) \right ]^{\frac{d-3}{2}} \\ \notag
f_{u} &\propto 2^{D-2} \left[ u-u^2 \right ]^{\frac{d-3}{2}} \\ \notag
f_{u} &\propto \left [ u-u^2 \right ]^{\frac{d-3}{2}} \\  \notag
f_{u} &\propto u^{\frac{d-1}{2} -1} (1-u)^{\frac{d-1}{2}-1}, \ given \frac{d-3}{2} = \frac{d-1}{2} -1  \\ \notag
f_{u} &\propto u^{\alpha -1} (1-u)^{\alpha-1}, with\  \alpha = \frac{d-1}{2} \notag
\end{align} \

The last line of equation (10) demonstrates that $u = (t + 1)/2$ has a symmetric beta distribution with  shape parameter $\alpha = (d-1)/2$. As $\alpha$ is a continuous variable  it is possible to speak about non integer values of intrinsic dimension d given $ \alpha = \frac{d-1}{2}$. \

The four references state that when d = 2 the beta density is an arcsine, d = 3 is an uniform, and when d = 4 the density function is a Wigner semi-circle.  It is apparent  that the extrinsic dimensions involved in the histogram model given in section (5) are much greater than d = 2,3, or 4. It will be shown, in the following section, that the intrinsic dimension of the rows of the symmetric circulant matrix match the values given by the four references. \

\section{Statistical Estimates  Derived from the Eigenanalysis of the Symmetric Circulant Matrix} \

A single autocorrelated long memory time series has a limited amount of statistical information. Sampling many independently produced, equal length, and standardized time series is one solution to derive statistical trends. The histograms of the probability density functions produced by the program in section 5.1 were generated from 200 samples of time series sequences. Another approach is to derive the ratio of variance to squared range of many standardized  time series, as demonstrated in section 6.1. Both variance and range are easily measured from the concatenation of many time series, and the ratio of variance to squared range determines the beta parameter $\alpha$, and  intrinsic dimension d.

Another solution is to estimate the statistical trends directly from the eigenvalues of the symmetric circulant matrix used to generate the time series. The  eigenvalues of the $rn \times rn$ matrix are indexed with respect to magnitude,  $\lambda_{1}$ being the largest, and $\lambda_{rn}$ being the smallest. What follows is a demonstration of the statistical trends that may be derived from the eigenvalues.
 
\subsection{Estimates of Intrinsic Dimension,  Beta Distribution Shape Parameter, and Principle Components } \

Ipsen and Saibaba \cite{iI24} define intrinsic dimension of a matrix A as:

\begin{equation}
\frac {trace(A)} { ||A||_{2}}
\end{equation} \ \

Where trace(A) is the sum of the main upper left to lower right diagonal elements of matrix A, or equally the sum of the eigenvalues of matrix A. The denominator of equation(11) is the spectral norm of matrix A,  the largest eigenvalue, $\lambda_{1}$. Equation (11) is a function of the  eigenvalues of matrix "circ".
 
The value of intrinsic dimension is context sensitive \cite{dS76}, and for this paper equation(11) is a good first approximation. 
A better estimate for intrinsic dimension d, between 2 to 4, and $\alpha$, between 0.5 and 1.5,  is  equation(12), an empirically derived two step estimate.

\begin{align}
E&= \frac {trace(circ)}{2^{1/2} ||circ||_{2}} +1\\  \notag
d &\cong E- \frac {E-3}{4.2}\\  \notag
\alpha &\cong \frac {d-1} {2} \notag
\end{align} \ \

\begin{figure}[h]
\centering \includegraphics [scale = .65]{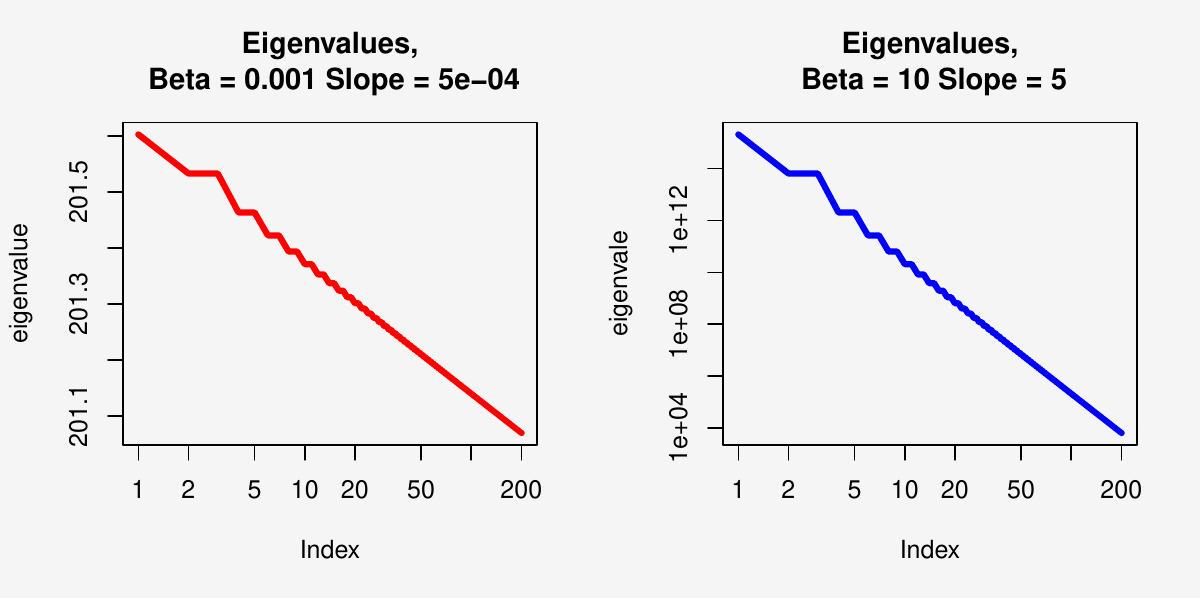}
\caption {Log Log Plots of Eigenvalues}
\label{Fi:003}
\end{figure}

Figure  3 are two log log plots of eigenvalues for a normal distribution beta = 0.001, and an arcsine distribution, beta = 10.  Replacing trace(A) in equation (11)  with the sum of eigenvalues, and $||A||_{2} $ by the largest eigenvalue $\lambda_{1}$ yields:

\begin{equation}
\frac { \lambda_1 + \sum_{i=2}^{rn}\lambda_i}{\lambda_1}
\end{equation}

Equation (13) while imprecise, is useful to observe trends at the extremes of the  range of "beta".
When beta approaches zero as in the left plot of figure 3, equation(13) approaches rn,  all Fourier moduli are approximately equal, and the distribution tends to  normal. When beta is ten, the right plot of figure 3, equation (13)  approaches unity, the sum of $\lambda_{1}$ and $\lambda_{2}$ dominate the sum of all eigenvalues, and they are the principle components of the time series. A time series with a beta of ten has an arcsine distribution,  the distribution of single frequency sinusoid. The bottom right plot of figure 1, with beta = 10,  is seen to be nearly sinusoidal. The plots in figure 3 represent the Fourier spectral trends of the long memory time series. The Fourier moduli of the time series are randomized by the product of "circ" with "epsilon", but on average, the spectral trend is an accurate statistic with a slope of $\beta/2$. The log log plot of the squared eigenvalues has a  slope equal to  $\beta$, the trend of the power spectral density.

\subsection{Eigenvalue  Estimate of The Variance of  Long Memory Time Series} \

Parseval's  theorem, when applied to a unitary discrete Fourier transform, with x the time domain variable, and X the Fourier frequency domain variable, can be expressed as:

\begin{equation}
\sum_{ j= 1 }^{n} |x_ j|^2 = \sum_{k = 1}^{n} |X_k|^2
\end{equation} \

In order to estimate the variance, it is necessary to remove the coefficient of the constant offset, zero frequency  Fourier term $X_1$.  Translating equation (14) to  an estimate of variance using the eigenvalues of "circ" that correspond to the Fourier coefficients $X_k$,  produces the following equation:

\begin{equation}
var(x) \cong \frac {1} {rn-1} \sum_{k = 2}^{rn}\lambda_{k} ^{2}
\end{equation} \

Given the estimates of $\alpha$ and variance, the range of the general beta distribution may be estimated. If the variance of "epsilon" is different than 1 then the right side of equation (15) needs to be multiplied by the new variance of "epsilon".

\subsection{The Condition Number of The Symmetric Matrix "Circ" and Geometric Collapse}  \  \

A  positive definite  matrix, a matrix with positive definite eigenvalues,  has a condition number $\kappa$ equal to the maximum eigenvalue divided by the minimum eigenvalue \cite{gS80}. 

\begin{equation}
\kappa = \frac {\lambda_{1}} {\lambda_{n}}
\end{equation} \

The matrix "circ", when "beta" = 0, has a distribution that  is normal, and  $\kappa$ = 1. When "beta" = 2.2,  the distribution is a Wigner semi-circle, and $\kappa = 3.4 \times 10^2$. When "beta" = 10,  the distribution is an arcsine, and $\kappa = 3.2 \times10^{11}$.

The reduction of intrinsic dimension, increase in serial correlation of the time series,  decreasing significance of smaller eigenvalues, and  approach to singularity of the matrix "circ", follows the increase of the condition number $\kappa$. This is called the geometric collapse of the matrix. The intrinsic dimension of the long memory  time series is the number of significant Fourier terms. Figure 1 illustrates the gradual reduction of high frequency spectral components as "beta" increases. 

\begin{table}
\begin{center}
\begin{tabular}  { | c | c | c | c | c | }
\hline
$\beta$ 	& Statistic 	&  Eigenvalue & Mean Measured  &  Measured	\\
        	&             	 &  Estimated		& Statistic               	& Coefficient of		\\ 
       	 &             	 &   Statistic        	 &                            &Variation			\\  \hline
2.2	 &	d           &    4.17                  &     4.22                &       0.22                \\   \hline
3.0       &       d           &      2.93                &      3.21               &       0.23                \\ \hline
10	  & 	d           &     2.05                 &     2.00                & 	0.03                \\ \hline
2.2       &    $\alpha$  &    1.58                  &     1.61                &      0.29                  \\ \hline
3.0       &    $\alpha$   &   0.97                  &     1.10                &      0.34                  \\  \hline
10        &    $\alpha$   &    0.52                 &      0.50               &       0.05                  \\ \hline 
2.2      &  variance      &   1.72 E 3           &       1.67 E 3          &      0.71                   \\ \hline 
3.0      &  variance       &   9.62 E 4         &        8.80 E 4          &      0.88                   \\ \hline
10	&  variance	    &   1.03 E 21         &       1.03 E 21        & 0.96                      \\ \hline

\end{tabular}
\linebreak
\linebreak
\caption {Statistics Estimated from Eigenvalues, and Measured from 500 Long Memory Sample Sequences of Length 200, using the R code in subsection 7.4}\label{ta : first}
\end{center}
\end{table}

\subsection{ The R language Program to Evaluate the Estimates of Intrinsic Dimension, Beta Shape Parameter, and Variance Derived from the Eigenvalues of "Circ", and Measured from Sample Variance and Range Squared }\

\begin{verbatim}
#Estimates of d, alpha, and variance derived from eigenvalues
#and measured from variance and range squared
#Load package "magic" to produce circulant matrix
beta = 2.2
n = 200
alpha = 0
d = 0 
v = 0

q = seq(from = -1/2, to = -1/n, by = 1/n)
r = 1/(2*n)
s = -rev(q)
frequency = c(q,r,s)
rn = length(frequency)
density = (abs(frequency))^(-beta/2)
b = Mod((rn^(-1/2))*fft(density, inverse = TRUE))
circ = (rn^(-1/2)) * circulant(b)
e = eigen(circ)
val = e$values

# ESTIMATION OF D AND ALPHA FROM EIGENVALUES
E = (sum(val[1:rn])/((2^(1/2))*val[1]))+1
d_est = E - (E-3)/4.2
alpha_est = (d_est - 1)/2

#ESTIMATION OF VARIANCE FROM EIGENVALUES
var_est = (1/(rn-1)) * sum(val[2:rn]^2)

#LOOP TO MEASURE ACTUAL MEAN AND COEFFICIENT
#OF VARIATION OF D, ALPHA, AND VARIANCE
for(i in 1:500) {
epsilon = rnorm(rn, mean = 0, sd = 1)
longmem = circ%*% epsilon

# 500 MEASURES OF ALPHA, D, AND VARIANCE
alpha[i]= (1/8) * ((max(longmem)-min(longmem))^
2 / var(longmem)) -1/2
d[i] = 2 * alpha[i] +1
v[i] = var(longmem) }

#ESTIMATES OF D, AND ALPHA FROM EIGENVALUES
alpha_est
d_est
formatC(var_est, format = "e")

#MEASURED MEAN, AND COEFFICIENT OF
#VARIATION FROM 500 SAMPLES
mean(alpha)
sd(alpha)/mean(alpha)
mean(d)
sd(d)/mean(d)
formatC(mean(v), format = "e")
sd(v)/mean(v)
\end{verbatim}

\section {Acknowledgment} \

The author found inspiration to pursue this topic and write this paper from the work of Herman Wold's 1938 book  "A Study In The Analysis Of Stationary Time Series" \cite{hW54}, as well as the ubiquity of 1/f noise in the outputs of astronomical instrumentation.

\pdfoutput =1

\vspace{.25in}

\end{document}